\input amstex 
\documentstyle{amsppt}
\input bull-ppt
\keyedby{bull304/mhm}

\topmatter
\cvol{27}
\cvolyear{1992}
\cmonth{October}
\cyear{1992}
\cvolno{2}
\cpgs{252-256}
\title Piercing convex sets \endtitle
\author Noga Alon and Daniel J. Kleitman \endauthor
\shortauthor{Noga Alon and D. J. Kleitman}
\address Department of Mathematics, Raymond and Beverly 
Sackler Faculty of
Exact 
Sciences, Tel Aviv University, Tel Aviv, Israel and 
Bellcore,
Morristown, New Jersey 07960\endaddress
\address Department of Mathematics, Massachusetts 
Institute of Technology, 
Cambridge, Massachusetts 02139\endaddress
\thanks Research supported in part by a United States 
Israel BSF Grant and by a
Bergmann Memorial Grant\endthanks
\thanks Supported in part by NSF grant DMS-9108403\endthanks
\date August 23, 1991\enddate
\abstract A family of sets has the $(p,q)$ property if 
among any $p$ members of
the family some $q$ have a nonempty intersection. It is 
shown that for every
$p\ge q\ge d+1$ there is a $c=c(p,q,d)<\infty$ such that 
for every family $\scr
F$ of compact, convex sets in $R^d$ that has the $(p,q)$ 
property there is a
set of at most $c$ points in $R^d$ that intersects each 
member of $\scr F$.
This extends Helly's Theorem and settles an old problem of 
Hadwiger and
Debrunner.\endabstract
\endtopmatter

\document

\heading 1. Introduction \endheading
For two integers $p\ge q$, a family of sets $\scr H$ has 
the $(p,q)$ {\it 
property\/} if among any $p$ members of the family some 
$q$ have a nonempty
intersection. $\scr H$ is $k$-{\it pierceable\/} if it can 
be split into $k$ or
fewer subfamilies, each having a nonempty intersection. 
The {\it piercing number\/}
of $\scr H$, denoted by $P(\scr H)$, is the minimum value 
of $k$ such that
$\scr H$ is $k$-pierceable. (If no such finite $k$ exists 
then $P(\scr
H)=\infty$.)

The classical theorem of Helly \cite {14} states that any 
family of compact
convex sets in $R^d$ that satisfies the $(d+1,d+
1)$-property is 1-pierceable.
Hadwiger and Debrunner considered the more general problem 
of studying the
piercing numbers of families $\scr F$ of compact, convex 
sets in $R^d$ that
satisfy the $(p,q)$ property. By considering the 
intersections of hyperplanes
in general position in $R^d$ with an appropriate box one 
easily checks that for
$q\le d$ the piercing number can be infinite, even if 
$p=q$. Thus we may assume
that $p\ge q\ge d+1$.

Let $M(p,q;d)$ denote the maximum possible piercing number 
(which is possibly
infinity) of a family of compact convex sets in $R^d$ with 
the
$(p,q)$-property. By Helly's Theorem, 
$$ M(d+1,d+1;d)=1$$
for all $d$ and, trivially, $M(p,q;d)\ge p-q+1$. Hadwiger 
and Debrunner \cite
{12} proved that for $p\ge q\ge d+1$ satisfying
$$p(d-1)<(q-1)d, \tag 1$$
this bound is tight, i.e., $M(p,q;d)=p-q+1$. In all other 
cases, it is not even known
if $M(p,q;d)$ is finite, and the question of deciding if 
this function  is
finite, raised by Hadwiger and Debrunner in 1957 in \cite 
{12}, remained open.
This question, which is usually referred to as the 
$(p,q)$-problem, is
considered in various survey articles and books, including 
\cite {13, 5, 8}.
The smallest case in which finiteness is unknown, which is 
pointed out in all
the above-mentioned articles, is the special case $p=4$, 
$q=3$, $d=2$. We note,
in all the cases where finiteness is known, that in fact 
$M(p,q;d)=p-q+1$ and
there are examples of Danzer and Gr\"unbaum (cf.\ \cite 
{13}) that show
$M(4,3;2)\ge 3>4-3+1$.

The $(p,q)$-problem received a considerable amount of 
attention, and finiteness
was proved for various restricted classes of convex sets, 
including the family
of parallelotopes with  edges parallel to the coordinate 
axes in $R^d$ \cite
{13, 19, 6}, families of homothetes of a convex set \cite 
{19}, and, using a
similar approach, families of convex sets with a certain 
``squareness'' 
property (\cite {9}, see also \cite {21}). 

Despite these efforts, the problem of deciding if 
$M(p,q;d)$ is finite remained
open for all values of $p\ge q\ge d+1$ that do not satisfy 
(1).

The purpose of the present note is to announce the 
solution of this problem,
established in the following 

\thm{Theorem 1.1} For every $p\ge q\ge d+1$ there is a 
$c=c(p,q,d)<\infty$ such
that $M(p,q;d)\le c${\rm ;} i.e., for every family $\scr 
F$ of compact, convex
sets in $R^d$ that has the $(p,q)$ property there is a set 
of at most $c$
points in $R^d$ that intersects each member of $\scr F$.
\ethm

The detailed proof will appear in \cite {3}. Here we 
briefly sketch the main
ideas. Three tools are applied: a fractional version of 
Helly's Theorem that
was first proved in \cite {15}, Farkas's Lemma (or Linear 
Programming Duality),
and a recent result proved in \cite {1}.

It may seem that there are almost no interesting families 
of compact convex
sets in $R^d$ that satisfy the  $(p,q)$-property for some 
$p\ge q\ge d+1$. A
large class of examples can be constructed as follows. Let 
$\mu$ be an 
arbitrary probability distribution on $R^d$, and let $\scr 
F$ be the family of
all compact convex sets $F$ in $R^d$ satisfying $\mu(F)\ge 
\varepsilon$. Since
the sum of the measures of any set of more than 
$d/\varepsilon$ such sets is
greater than $d$, it follows that if $p$ is the smallest 
integer strictly
larger than $d/\varepsilon$ then $\scr F$ has the $(p,d+
1)$ property. It then
follows that $P(\scr F)\le M(p,d+1;d+1)$, i.e., for every 
probability measure
in $R^d$ there is a set $X$ of at most $M(p,d+1;d+1 )$ 
points such that any
compact convex set in $R^d$ whose measure exceeds 
$\varepsilon$ intersects $X$.

The following theorem is an immediate consequence of 
Theorem 1.1.

\thm{Theorem 1.2} Let $\scr F$ be a family of compact 
convex sets in $R^d$, and
suppose that for every subfamily $\scr F'$ of cardinality 
$x$ of $\scr F$ the 
inequality $P(\scr F')<\lceil x/d\rceil $ holds\/{\rm ;} 
i.e., $\scr F'$ can be
pierced by less than $x/d$ points. Then $P(\scr F)\le 
M(x,d+1;d+1)$.
\ethm

Observe that in order to deduce a finite upper bound for 
the piercing number of
$\scr F$, the assumption that $\scr P(\scr F')<\lceil 
x/d\rceil$ cannot be
replaced by $P(\scr F')\le \lceil x/d\rceil$ as shown by 
an infinite family of
hyperplanes in general position (intersected with an 
appropriate box), whose
piercing number is infinite.

\heading 2. A sketch of the proofs \endheading
Since we do not try to optimize the constants here, and 
since obviously
$M(p,q;d)\le M(p,d+1;d)$ for all $p\ge q\ge d+1$, it 
suffices to prove an upper
bound for $M(p,d+1;d)$. Another simple observation is that 
by compactness we
can restrict our attention to finite families of convex 
sets.

Let $\scr F$ be a family of $n$ convex sets in $R^d$, and 
suppose that $\scr F$
has the $(p,d+1)$ property. Our objective is to find an 
upper bound for the
piercing number $P(\scr F)$ of $\scr F$, where the bound 
depends only on $p$
and $d$. It is convenient to describe the ideas in three 
subsections.

\subheading{\num{2.1.} A fractional version of Helly's 
Theorem} 
Katchalski and Liu \cite {15} proved the following result, 
which can be viewed
as a fractional version of Helly's Theorem.

\thm{Theorem 2.1 \cite {15}\rm} For every $0<\alpha\le 1$ 
and every $d$ there is a
$\delta=\delta(\alpha,d)>0$ such that for every $n\ge d+
1$, every family of $n$
convex sets in $R^d$ that contains at least $\alpha\binom 
{n}{d+1}$
intersecting subfamilies of cardinality $d+1$ contains an 
intersecting
subfamily of at least $\delta n$ of the sets.
\ethm

Notice that Helly's Theorem is equivalent to the statement 
that in the above 
theorem $\delta(1,d)=1$. 

A sharp  quantitative version of this theorem was proved 
by Kalai \cite {16}
and, independently, by Eckhoff \cite {7}. See also \cite 
{2} for a very short
proof. All these proofs rely on Wegner's Theorem \cite 
{20} that asserts that
the nerve of a family of convex sets in $R^d$ is 
$d$-collapsible.

Theorem 2.1, together with a simple probabilistic 
argument,  can be applied to
prove 

\thm{Lemma 2.2} For every $p\ge d+1$ there is a positive 
constant
$\beta=\beta(p,d)$ with the following property. Let $\scr 
F=\{A_1,\dotsc,
A_n\}$ be a family of $n$ convex sets in $R^d$ that has 
the $(p,d+1)$ property.
Let $a_i$ be nonnegative integers, define 
$m=\sum^n_{i=1}a_i$, and let $\scr G$
be the family of cardinality $m$ consisting of $a_i$ 
copies of $A_i$ for $1\le
i\le n$. Then there is a point $x$ in $R^d$ that belongs 
to at least $\beta m$
members of $\scr G$.
\ethm

\subheading{\num{2.2.} Farkas's Lemma and a lemma on 
hypergraphs}
The following is a known variant of the well-known lemma 
of Farkas (cf.\ \cite
{17, p.\ 90}).

\thm{Lemma 2.3} Let $A$ be a real matrix and $b$ a real 
{\rm (}column\/{\rm )}
vector. Then the system $Ax\le b$ has a solution $x\ge 0$ 
if and only if for
every {\rm (}row\/{\rm )} vector $y\ge 0$ that satisfies 
$yA\ge 0$ the 
inequality $yb\ge 0$ holds.
\ethm

This lemma can be used to prove the following.

\thm{Corollary 2.4} Let $H=(V,E)$ be a hypergraph and let 
$0\le \gamma\le 1$ be
a real. Then the following two conditions are 
equivalent\/{\rm :}
\roster
\item "(i)" There exists a weight function $f\:V\mapsto R^+
$ satisfying
$\sum_{v\in V}f(v)=1$ and $\sum_{v\in e}f(v)\ge \gamma$ 
for all $e\in E$.
\item "(ii)" For every function $g\:E\mapsto R^+$ there is 
a vertex $v\in V$
such that $\sum_{e;v\in e}g(e)\ge \gamma \sum_{e\in E}g(e)$.
\endroster
\ethm

By the last corollary and Lemma 2.2 one can prove the 
following result.

\thm{Corollary 2.5} Suppose $p\ge d+1$ and let 
$\beta=\beta(p,d)$ be the
constant from Lemma {\rm 2.2}. Then for every family $\scr 
F=\{A_1,\dotsc, 
A_n\}$ of $n$ convex sets in $R^d$ with the $(p,d+1)$ 
property there is a
finite {\rm (}multi\/{\rm -)}set $Y\subset  R^d$ such that 
$|Y\cap A_i|\ge 
\beta|Y|$ for all $1\le i\le n$.
\ethm

\subheading{\num{2.3.} Weak $\varepsilon$-nets for convex 
sets}
The following result is proved in \cite {1}.

\thm{Theorem 2.6 \cite {1}\rm} For every real 
$0<\varepsilon<1$ and every integer
$d$ there exists a constant  $b=b(\varepsilon,d)$ such 
that the following
holds\/{\rm :}

For every $m$ and every multiset $Y$ of $m$ points in 
$R^d$, there is a subset
$X$ of at most $b$ points in $R^d$ such that the convex 
hull of any subset of
$\varepsilon m$ members of $Y$ contains at least one point 
of $X$.
\ethm

Several arguments that supply various upper bounds for 
$b(\varepsilon,d)$ are
given in \cite {1}. The simplest one is based on a result 
of B\'ar\'any 
\cite {4} whose proof is based on a deep result of 
Tverberg \cite {8}.

Theorem 1.1 follows from the above results quite easily. 
Let $\scr
F=\{A_1,\dotsc, 
A_n\}$ be a family of $n$ convex sets in $R^d$ with the
$(p,d+1)$ property, where $p\ge d+1$. By Corollary 2.5 
there is a finite
(multi-)set $Y\subset  R^d$ such that $|Y\cap A_i|\ge 
\beta|Y|$ for all $1\le
i\le n$, where  $\beta=\beta(p,d)$ is as in Lemma 2.2. By  
Theorem 2.6 there is
a set $X$ of at most $b(\beta,d)$ points in $R^d$ such 
that the convex hull of
any set of $\beta|Y|$ members of $Y$ contains at least one 
point of $X$. Since
each member of $\scr F$ contains at least $\beta|Y|$ 
points in $Y$, it must
contain at least one point of $X$. Therefore, $P(\scr 
F)\le |X|\le
b(\beta(p,d),d)$, completing the proof.

The detailed proofs of the above lemmas and corollaries, 
as well as some
methods to improve the estimates for the numbers 
$M(p,q;d)$, will appear in
\cite {3}.

\enddocument